\documentclass[11pt,reqno]{amsart}
\usepackage{amssymb}
\usepackage{mathptmx}
\usepackage{mathtools}
\usepackage{mathrsfs}
\usepackage[export]{adjustbox}
\usepackage[all]{xy}
\usepackage{stmaryrd}
\usepackage{bbm}
\usepackage{fancyhdr}
\usepackage{xcolor}
\usepackage{hyperref}
\usepackage{tikz-cd}
\usepackage{cite}
\usepackage{amsmath,amsfonts}
\usepackage{graphicx}
\usepackage{wrapfig}
\usepackage{array}
\usepackage{amsthm}
\usepackage[left=1.2in, right=1.2in, top=1in, bottom=1in]{geometry}
\usepackage{etoolbox}

\DeclareMathOperator{\piop}{\pi}

\patchcmd{\section}{\scshape}{\bfseries}{}{}
\makeatletter
\renewcommand{\@secnumfont}{\bfseries}
\makeatother
\xyoption{all}
\newcommand{\mycomment}[1]{}

\theoremstyle{definition}
\newtheorem{mydef}{\textbf{Definition}}[section]
\newtheorem{myeg}[mydef]{\textbf{Example}}

\theoremstyle{plain}

\newtheorem{mythm}[mydef]{\textbf{Theorem}}
\newtheorem{lem}[mydef]{\textbf{Lemma}}

\begin{document}

\title{Blocks in Finite Hyperfields}

	\author{David Hobby}
	\address{State University of New York at New Paltz, NY, USA}
	\curraddr{}
	\email{hobbyd@newpaltz.edu}

	
	\makeatletter
	\@namedef{subjclassname@2020}{%
		\textup{2020} Mathematics Subject Classification}
	\makeatother
	
	\date{}
	
	\dedicatory{}

\begin{abstract}	
This paper studies the structure of finite hyperfields $H$, and finds a subtle pattern in their addition operation.  
Consider the class $\mathcal{H}$ of all hyperfields with a given multiplicative group on $H^\times = H - \{0\}$ and given value of $-1$.
Then the addition of hyperfields in this class is determined by the set of pairs $(x,y)$ with $y \in x+1$ for $x,y \in H^\times$.
There are {\em blocks} of such pairs, where $(x_0,y_0)$ and $(x_1,y_1)$ are in the same block iff any hyperfield with $y_0 \in x_0 + 1$ also has $y_1 \in x_1 + 1$.
The theory of these blocks is developed, they can easily be computed without using hyperfields.
Exploiting this theory of blocks would greatly speed up future computer searches for small hyperfields.  
The theory of blocks is then used to show that the number of nonquotient hyperfields of size $n$ grows exponentially with $n$, and that for even $n$ most hyperfields are nonquotient.
A final application shows that a large class of finite hyperfields has the FETVINS property, meaning that systems of linear homogeneous equations with fewer equations than variables always have nontrivial solutions.
\end{abstract}

\maketitle


\section{Introduction}

Hyperfields are algebraic structures which assume the same axioms as fields except that one allows addition to be \emph{multi-valued}. A multi-valued addition may seem exotic, but the construction is natural and useful.  For any given field $F$ and multiplicative subgroup $G\leq F^\times$, the pair $(F,G)$ naturally defines a hyperfield $F/G$, called a quotient hyperfield (see Example \ref{quotient example}). In fact, Krasner \cite{krasner1956approximation} utilized the notion of quotient hyperfields to define a notion of limits of local fields, and investigated their algebraic extensions and Galois groups in positive characteristic. In doing so, Krasner \cite{krasner1983class } noticed that in order to show that the theory of hyperfields (or hyperrings) could deliver new results, one needed to find examples of non-quotient hyperfields. Shortly after that, C.G.~Massouros \cite{Mas85a, massouros1985theory} and A.~Nakassis \cite{nakassis1988recent} found examples of non-quotient hyperfields. Over time, more hyperfields were proved to be non-quotient.  We refer the reader to \cite{massouros2023on} for an overview.

Hyperfields are being actively studied in connection to many areas of mathematics. For instance, several authors employ hyperfields in number theory and algebraic geometry: A.~Connes and C.~Consani \cite{con3}, O.~Viro \cite{viro2010}, O.~Lorscheid \cite{lorscheid2022tropical}, and J.~Jun \cite{jun2017geometry}. 

The flexibility of the structure of hyperfields allows one to unify various mathematical objects under the same framework. In fact, in their celebrated work \cite{baker2017matroids}, M.~Baker and N.~Bowler unify various generalizations of matroids by using hyperfields, thereby allowing one to prove theorems concerning all these generalizations of matroids at the same time. 

With the same motivation, in \cite{baker2021descartes}, M.~Baker and  O.~Lorscheid  develop a theory of multiplicities of roots for polynomials over hyperfields to unify both Descartes' rule of signs and Newton's polygon rule. This naturally leads one to the work of M.~Baker and T.~Zhang \cite{baker2022some}, where they unify and provide conceptual proof for various results about ranks of matrices that are scattered throughout the literature by using the framework of matroids over hyperfields. 

Section 2 reviews the basic theory of hyperfields, covering the axioms, the quotient construction, and standard examples.

Section 3 develops the theory of blocks.
To make a hyperfield $H$, we may fix an abelian multiplicative group, which will give multiplication on $H^\times = H - \{0\}$.  
The multiplication of $H$ is then determined, since $x \cdot 0 = 0 \cdot x = 0$ follows from the axioms.
Once we know which element of $H$ is $-1$, all additions involving $0$ are determined, for the axioms give us 
$x + 0 = 0 + x = x$, and also that $0$ is in the set $x+y$ iff $y = (-1) \cdot x$.

This leaves addition on $H^\times$ to be determined.
The sum $x + y$ is equal to $y \cdot (y^{-1}x + 1)$, so it can be obtained using distributivity from the values of $x+1$ for $x \in H^\times$.
And finding the result of adding $1$ to $x \in H^\times$ is equivalent to knowing which pairs $(x,y)$ are in the relation $\piop$, which we define by   
$(x,y) \in \piop$ iff $y \in x+1$.

This is a complicated way to represent addition, but the payoff is that pairs in $H^\times \times H^\times$ come in {\em blocks}, where two pairs $(x_0,y_0)$ and $(x_1,y_1)$ are in the same block iff any hyperfield with $(x_0,y_0)$ in its $\pi$ relation also has $(x_1,y_1)$ in its $\pi$.

The theory of these blocks is developed, they can easily be computed without using hyperfields.
Since blocks may have up to $6$ elements, this promises a way to easily find all hyperfields of size $n$ for values of $n$ that were previously inaccessible.
For example, an $11$-element hyperfield would have a multiplicative group isomorphic to ${\mathbb Z}_{10}$.
It turns out that with both possible values of $-1$, the $\piop$ relation has 22 blocks, giving us $2^{22}$ possible choices of blocks to include and around $4$ million candidate hyperfields to examine.
This can be easily done by any computer.
The current best exhaustive listing of hyperfields is for those with $7$ elements and multiplicative group isomorphic to  ${\mathbb Z}_{6}$, in \cite{MassourosOrder7}.

It is not true that every structure with a $\piop$ relation that is a union of blocks is a hyperfield.
But it is true in the case where the structure is what we call {\em ample}.
This result is our Theorem \ref{large gives hyperfield theorem}, and ample hyperfields are defined in Definition \ref{ample definition}.
Crudely stated, a hyperfield $H$ is ample if at least half of the pairs in $H^\times \times H^\times$ are in its $\piop$ relation.

Section 4 uses the theory of blocks to show that the number of nonquotient hyperfields of size $n$ grows exponentially with $n$, and that for even $n$ most hyperfields are nonquotient.
This is a partial answer to Question 5.4 of \cite{BakerJin2021}.
In that paper, M. Baker and T. Jin let $H_n$ and $Q_n$ denote the number of isomorphism classes of $n$-element hyperfields and $n$-element quotient hyperfields, respectively.
They ask about the limit as n goes to infinity of the quantity $Q_n / H_n$, and conjecture the limit is zero.
We show in Theorem \ref{limit of ratio theorem} that the limit is zero when $n$ is restricted to even values.

Section 5 shows that all ample hyperfields have FETVINS.

In Question 5.13 of \cite{baker2022some}, M. Baker and T. Zhang ask if a system of homogeneous linear equations over a hyperfield $H$ with more unknowns than equations always has a nontrivial solution.
This is a very natural question, which is seemingly difficult to answer.
In the paper \cite{hobbyjun2024}, I and J. Jun made some progress on this question.  
It is useful to have a name for the property that these systems always have nontrivial solutions.
We used the acronym FETVINS, for ``Fewer Equations Than Variables Implies Nontrivial Solutions''.

As Baker and Zhang observe, all quotient hyperfields have FETVINS, they inherit it from the original fields.
So it was natural to investigate the question for the remaining hyperfields, the nonquotient ones.
In \cite{hobbyjun2024}, we showed that most of the original examples of nonquotient hyperfields had FETVINS.
This used the fact that they had a property similar to being ample.
Now that the property of ``ampleness'' has been discovered , it is natural to try to show that ample hyperfields have FETVINS.
Indeed they do, this is Theorem \ref{ample implies FETVINS theorem}.
The proof is similar to the analogous one in \cite{hobbyjun2024}, but fortunately simpler.

\section{Preliminaries}
It is convenient to sometimes write $S \ni x$ for $x \in S$, so we will do so.

\begin{mydef} \label{hyperaddition definition}
Let $H$ be a nonempty set. By a hyperaddition, we mean a function
\begin{equation}
+:H \times H \to \mathcal{P}^*(H), 
\end{equation}
where $\mathcal{P}^*(H)$ is the set of nonempty subsets of $H$, satisfying the following conditions:
\begin{enumerate}
    \item
$a+b=b+a$. 
\item 
$(a+b)+c=a+(b+c)$\footnote{For any subset $A\subseteq H$ and an element $x \in H$, we let $A+x=\cup_{a \in A}(a+x)$.}
\end{enumerate}
\end{mydef}

When $a+b=\{c\}$, we simply write $a+b=c$. 

\begin{mydef} \label{hypergroup definition}
By a hypergroup, we mean a set $H$ with  hyperaddition which satisfies the following conditions:
\begin{enumerate}
    \item (Zero) \quad
$\exists~!~ 0\in H$ such that $0+h=h$ for all $h \in H$
\item  (Negative) \quad
For each $x \in H$, $\exists~!~y \in H$ such that $0 \in x+y$. We denote $y:=-x$. 
\item  (Reversibility) \quad
If $x \in y+z$, then $z \in (x-y):=x+(-y)$. 
\end{enumerate}
\end{mydef}

\begin{mydef} \label{hyperring definition}
By a hyperring, we mean a set $R$ with hyperaddition $+$ and usual multiplication $\cdot$ satisfying the following:
\begin{enumerate}
    \item 
$(H,+)$ is a hypergroup. 
\item 
$(H,\cdot)$ is a multiplicative monoid. 
\item 
$a(b+c)=ab+ac$ for all $a,b,c \in R$. 
\item 
$r\cdot 0 = 0$ for all $r \in R$. 
\end{enumerate}
If $H - \{0\}$ is a multiplicative group, then $H$ is said to be a hyperfield. 
\end{mydef}

In the presence of the other hyperfield axioms, the Reversibility Axiom is redundant.  
This has been noted before by various authors, such as C. G. Massouros and G. G. Massouros in \cite{MassourosOrder7}, and 
J. Maxwell and B. Smith in \cite{ConvexGeometry}.
The situation is complicated because the axioms for hyperfields vary a bit between authors.
We will need this fact later, and present it here as a lemma.

\begin{lem} \label{do not need reversibility lemma}
The reversibility axiom, (3) of Definition \ref{hypergroup definition}, is implied by the other hyperfield axioms.
\end{lem}

\begin{proof}
Given a structure $H$ that satisfies all the other hyperfield axioms, let $x, y, z \in H$ be given.
We must show $x \in y+z$ iff $x-y \ni z$.

Assuming $x \in y+z$, we have $0 \in x-x \subseteq (y+z) - x = z + (-x + y)$.
Thus $-z \in -x+y$, so $-x + y \ni -z$ and 
$x-y \ni z$.
The other direction is similar, given 
$x-y \ni z$, we have 
$0 \in x-y-z$ so $-x \in -y-z$ and $x \in y+z$.
\end{proof}

\begin{myeg}[Krasner hyperfield]
Let $\mathbb{K}=\{0,1\}$. We define the following addition:
\[
0+1=1, \quad 0+0=0, \quad 1+1=\{0,1\}
\]
Multiplication is the same as in $\mathbb{F}_2$. Then $\mathbb{K}$ is a hyperfield, called the Krasner hyperfield. 
\end{myeg}

\begin{myeg}[Sign hyperfield]
Let $\mathbb{S}=\{-1,0,1\}$. We define the following addition:
\[
0+1=1, \quad 0+(-1)=-1,\quad 0+0=0, \quad 1+1=1, \quad (-1)+(-1)=-1, \quad 1+(-1)=\mathbb{S}. 
\]
Multiplication is that same as in $\mathbb{F}_3$. Then $\mathbb{S}$ is a hyperfield, called the sign hyperfield. 
\end{myeg}

\begin{myeg}[Tropical hyperfield]
Let $\mathbb{T}:=\mathbb{R}\cup \{-\infty\}$. Hyperaddition for $\mathbb{T}$ is as follows: for $a,b \in \mathbb{T}$, 
\[
a+b=\begin{cases}
\max\{a,b\} \textrm{ if $a\neq b$,} \\
[-\infty, a] \textrm{  if $a=b$}
\end{cases}
\]
Multiplication is the usual addition of real numbers with $a\cdot(-\infty)=-\infty$. 
\end{myeg}

\begin{myeg}[Phase hyperfield]
Let $\mathbb{P}=S^1\cup \{0\}$, where $S^1$ is the complex unit circle. Multiplication is usual multiplication of complex numbers, and addition is defined as follows: for $a, b \in S^1$
\begin{equation}
a+b=\begin{cases}
\{0,a,-a\} \textrm{ if $b=-a$ }, \\
\textrm{all points in the shorter of the two open arcs of $S^1$ connecting $a$ and $b$ if $b\neq a$.}
\end{cases}
\end{equation}
Also $a+0=a$ for all $a \in \mathbb{P}$. 
\end{myeg}

\begin{myeg} \label{quotient example}
In general, let $A$ be a commutative ring and $G \leq A^\times$, the multiplicative group of the units of $A$. Then $G$ naturally acts on $A$ by multiplication and we obtain the set $A/G$ of equivalence classes. Let $[a]$ be the equivalence class of $a \in A$. One defines the following hyperaddition:
\[
[a]+[b]:=\{[c] \mid c=ag_1+bg_2 \textrm{ for some } g_1,g_2 \in G\}
\]
Multiplication is $[a]\cdot[b]=[ab]$. With these two operations, $A/G$ becomes a hyperring and if $A$ is a field, then $A/G$ is a hyperfield. 
\end{myeg}

The Phase hyperfield is also a nice example of this quotient construction.  
We let $F$ be the field of complex numbers, and let $G$ be the set of positive reals. 
Then $G$ is a subgroup of $F^\times = {\mathbb C} - \{0\}$, and the cosets of $G$ are sets of the form 
$G_\theta = \{r e^{i \theta} \colon r \in {\mathbb R}, r > 0\}$
for $0 \leq \theta < 2 \pi$.
The product of two cosets is 
$G_\alpha \cdot G_\beta =  \{ r e^{i \alpha} + r' e^{i \beta} \colon r,r' \in {\mathbb R}, r,r' > 0\} =
\{ s e^{i (\alpha + \beta)} \colon s \in {\mathbb R}, s > 0\}$
The sum of two cosets is 
$G_\alpha + G_\beta = \bigcup \{ r e^{i \alpha} + r' e^{i \beta} \colon r,r' \in {\mathbb R}, r,r' > 0\} $,
viewed as a set of cosets.
This quotient hyperfield is isomorphic to the Phase hyperfield.

\section{The theory of blocks}
As part of a research project with students, the author looked over  
the paper \cite{Ameri2020} of M. Eyvazi, R. Ameri and S. Hoskova-Mayerova.
It had an exhaustive listing of small hyperfields found by computer search.
He was struck by the fact that in many of these hyperfields sums of two elements were relatively large sets.

The somewhat cynical response to this is the observation that to show associativity, it suffices to have all sums of three elements equal to $H$.
Slightly modifying this, we are led to the following lemma.

\begin{lem} \label{Large sum lemma}
Let $H$ be a finite structure satisfying all the hyperfield properties except possibly associativity of addition.
Assume that $m$ is the largest number so that for all $w \in H^\times$ we have $|(w+1)-\{0\}| \geq m$.
Also let $k$ be the largest number so that for all $w \in H^\times$ we have that $|\{u \in H^\times \colon w \in (u+1)\}|$ is at least $k$.
Then if $m+k > |H^\times|$, we have for all $x,y,z \in H^\times$, that
$H^\times$ is a subset of both $(x+y)+z$ and $x+(y+z)$.
\end{lem}

\begin{proof}
Let $H$ be such a structure, where $m+k > |H^\times|$.
Since $H$ is commutative and $x,y,z \in H^\times$ are arbitrary, it suffices to show $H^\times \subseteq (x+y)+z$ for any given $x$, $y$ and $z$.
Indeed, since 
$|(x+y)+z| = |(y^{-1}x + 1) + y^{-1}z|$ by distributivity, we may assume $y = 1$.
So given $v \in H^\times$, we will show 
$v \in (x+1)+z$.
We define $K_v \subseteq H^\times$ to be
$\{w \in H^\times \colon v \in w+z\}$.
Now we have $v \in w+z$ iff 
$vz^{-1} \in wz^{-1} + 1$, so 
$|K_v| = |\{u \colon vz^{-1} \in u+1 \}| \geq k$.
Then there is some $s \in H^\times$ that is in both 
$x+1$ and $K_v$, since these sets have sizes $m$ and $k$, where $m+k > |H^\times|$.
But then $v \in s + z \subseteq (x+1)+z$, as desired.
\end{proof}

Under the assumptions of Lemma \ref{Large sum lemma}, we have for all $x,y,z \in H^\times$, that the sets $(x+y)+z$ and $x+(y+z)$ agree on all of $H$ except $0$.
So it is natural to investigate when we have $0 \in (x+y)+z$ iff $0 \in x+(y+z)$ for all $x,y,z \in H^\times$.
We have that $0 \in (x+y)+z$ iff
$0 \in (y^{-1}x + 1)+y^{-1}z$ and
$0 \in x+(y+z)$ iff
$0 \in y^{-1}x + (1+y^{-1}z)$.
Replacing $y^{-1}x$ with $x$ and $y^{-1}z$ with $z$, we have that addition on $H^\times$ is associative just in case we have for all $x,z \in H^\times$ that
$0 \in (x+1)+z$ iff $0 \in x+(1+z)$.

Now $0 \in (x+1)+z$ iff $-z \in (x+1)$, and
$0 \in x+(1+z)$ iff $-x \in (1+z)$.
Thus we need $-z \in (x+1)$ iff $-x \in (1+z)$ for all $x,z \in H^\times$.
This is better expressed as a closure condition on the set of pairs $(x,z)$ with
$z \in x+1$.  

As has been noted many times, for example in \cite{Ameri2020}, the addition of a hyperfield is determined by the results of adding $1$.
So it is natural to ask when a ``plus one contains'' relation $\piop$ defines a hyperfield.

\begin{mydef} \label{H_pi definition}
Let $H^\times$ be an abelian group, written multiplicatively, and let $-1$ be a distinguished element of $H^\times$ that has order $1$ or $2$.
Given a relation $\piop$ on $H^\times$, we define the structure $H_\pi$ as follows.
The universe of $H_\pi$ is $H^\times$ together with a new element $0$.
Multiplication on $H_\pi$ is defined to be that of $H^\times$ with the added condition that $0\cdot x = x \cdot 0 = 0$ for all $x \in H_\pi$.
Addition with $0$ is defined by $0 + x = x + 0 = x$ for all $x \in H_\pi$.
If $x$ and $y$ are both nonzero, we define 
$x+y$ to be $y$ times the set $P(y^{-1}x)$, where $P(z)$ is $\{ w \colon z \piop w\}$ for $z \neq -1$ and   
$P(-1) = \{ w \colon -1 \piop w\} \cup \{0\}$.
Thus addition is determined by $\piop$ once $-1$ is known.
\end{mydef}

Some conditions on $\piop$ are obvious. 
Since no sums in a hyperfield can yield the empty set, we must have that for all $x \in H^\times$ the set 
$P(x) \supseteq \{ y \colon x \piop y \}$ defined in the previous paragraph must be nonempty.
Another condition is needed so that addition is commutative.

\begin{mydef}
A relation $\piop$ on $H^\times$ is {\em consistent} iff whenever 
$x \in H^\times - \{1\}$
we have $x \piop y$ iff $x^{-1} \piop x^{-1}y$.
\end{mydef}

If addition is commutative in $H_\pi$, then the relation $\piop$ must be consistent.
For $x \piop y$ implies $y \in x + 1$, and multiplying by $x^{-1}$ gives 
$x^{-1}y \in 1 + x^{-1} = x^{-1} + 1$, 
so $x^{-1} \piop x^{-1}y$.
We will later show that a consistent $\piop$ gives a commutative addition.

Recall that we also had a condition that was required for addition to be associative.
This came from the requirement that 
for all $x,z \in H^\times$ that
$0 \in (x+1)+z$ iff $0 \in x+(1+z)$, which reduced to 
$-z \in (x+1)$ iff $-x \in (1+z)$.
Expressing this in terms of $\piop$, we need
$x \piop -z$ iff $-x \piop z$.
In other words, $\piop$ is closed under 
{\em reversal/negation}.
This discussion gives us the following theorem.

\begin{mythm}  \label{Closure theorem}
The $\piop$ relation of any hyperfield must be consistent and closed under reversal/negation.
\end{mythm}

This is a powerful theorem, because it greatly speeds up the search for finite hyperfields.
Consistency is also a closure condition, so the combined requirements of consistency and closure under reversal/negation allow us to replace single locations in the incidence matrix of $\piop$ with groups of several locations that all have the same value.
An example is in order.

We consider $4$-element hyperfields, so $H^\times$ is the $3$-element group 
$\{1,a,a^2\}$, and $-1$ is $1$, the sole element of order $1$ or $2$.
So $0$ is in $1+1$, but that is not part of the definition of $\piop$.
Suppose that $1 \piop 1$.
Then consistency gives us nothing, since it only applies to $x \piop y$ where $x \neq 1$.
And reversal/negation gives us 
$-1 \piop -1$, which is also $1 \piop 1$.
Thus the pair $(1,1)$ in $\piop$ has the set 
$A = \{(1,1)\}$ as its closure.

Next consider $1 \piop a$.
Consistency gives us $1 \piop a$ again, but 
reversal/negation gives $a \piop 1$.
Applying reversal/negation to $a \piop 1$ is pointless, it takes us back to $1 \piop a$.
But using consistency turns $a \piop 1$
into $a^{-1} \piop a^{-1} 1$, or 
$a^2 \piop a^2$.
Reversal/negation of $a^2 \piop a^2$ just gives us $a^2 \piop a^2$, so we have a second 
closed set of pairs,
$B = \{(1,a),(a,1),(a^2,a^2)\}$.

So we go to the next pair that is not already in a closed set, and consider 
$(1,a^2)$ or $1 \piop a^2$.
Reversal/negation gives $a^2 \piop 1$, and consistency applied to that gives 
$a^{-2} \piop a^{-2} 1$ or $a \piop a$.
The reader may check that no other pairs are produced, and we have the closed set 
$C = \{(1,a^2),(a^2,1),(a,a)\}$.

Two pairs remain, $(a,a^2)$ and $(a^2,a)$.
Reversal/negation switches them, as does consistency.
So we have a final closed set of pairs,
$D = \{(a,a^2),(a^2,a)\}$.

\begin{mydef} \label{block definition}
Given an abelian group $H^\times$ with distinguished element $-1$, {\em blocks} are sets of pairs of elements of $H^\times$ that are closed under consistency and reversal/negation. 
\end{mydef}

We can nicely represent the blocks for the group $H^\times$ in the following table.
\begin{center}
$\begin{array}{|c||c|c|c| } 
\hline
\piop & 1 & a & a^2 \\
\hline
\hline
1 & A & B & C \\
\hline
a & B & C & D \\
\hline
a^2 & C & D & B \\
\hline
\end{array}$
\end{center}

Every $4$-element hyperfield corresponds to a collection of blocks, although not all collections of blocks correspond to hyperfields.
We will list collections of blocks by sequences of capital letters.
The empty sequence, and also the sequence A do not give hyperfields, since some additions would be the empty set.
The sequence D gives the $4$-element field with set of elements $\{0,1,a,a^2\}$
The addition is determined as follows.
$1+1$ would be the empty set, except $-1 = 1$ so $1+1 = \{0\}$.
$a+1$ is $\{a^2\}$, or expressed more normally, $a + 1 = a^2$.
Similarly, $a^2 +1 = a$.
Together with distributivity, these facts determine the hyperfield.

The sequence BD gives us that $1+1 = \{a\}\cup\{0\} = \{0,a\}$, that $a+1 = \{1,a^2\}$, and that $a^2 +1 = \{a,a^2\}$.
This is isomorphic to the quotient hyperfield $\mathbb{Z}_7 / \langle 6 \rangle$, which is called $HF_{42}$ in
\cite{Ameri2020} and covered as (2) in the enumeration of hyperfields of order $4$ in 
\cite{BakerJin2021}.
The sequence CD gives us a hyperfield isomorphic to the one given by BD.
This is easily seen since $0$ and $1$ are determined but switching $a$ and $a^2$ gives an automorphism of the multiplicative group.
The action of this automorphism on $\piop$ fixes the block A, and swaps the pairs 
$(a,a^2)$ and $(a^2,a)$ leaving the block D invariant.
It also swaps the pairs $(1,a^2)$ and $(1,a)$, the pairs $(a,a)$ and $(a^2,a^2)$, and the pairs $(a^2,1)$ and $(a,1)$.
This results in exchanging the blocks C and B, showing that the hyperfields given by BD and CD are isomorphic.

Following the enumeration in Case 3 of \cite{BakerJin2021}, the remaining $4$-element hyperfields are given by the following blocks.
The sequence BCD gives the hyperfield listed as (4), and it is isomorphic to the quotient 
$\mathbb{Z}_{13} / \langle 8 \rangle$.
The sequence ABCD gives the hyperfield listed as (5), isomorphic to the quotient 
$\mathbb{Z}_{19} / \langle 8 \rangle$.
The sequence BC gives the hyperfield listed as (6), which is nonquotient.
The sequences ABD and ACD give the two isomorphic hyperfields listed as (7) and (8), which are nonquotient.
The sequence ABC gives the nonquotient hyperfield listed as (9).

The remaining sequences are B, C, AB, AC, BC and AD.  
They fail to give hyperfields since addition is not associative in the resulting structures.
As we will see, this is typical.
Sequences with many blocks give hyperfields, while most ``sparse'' sequences do not.

\begin{mythm} \label{block size at most 6 Theorem}
Every block that is closed under consistency and reversal/negation consists of $6$ or fewer pairs.
And every block that has a pair of the form $(1,y)$ in it contains at most $3$ pairs, $(1,y)$, $(-y,-1)$ and $(-y^{-1},y^{-1})$.
\end{mythm}

\begin{proof}
Consider an arbitrary initial pair $(x,y)$ in a block.
Doing reversal/negation to $(x,y)$ gives $(-y,-x)$, and doing it again returns us to $(x,y)$.
Similarly, when $x \neq 1$ and consistency can be applied, two applications return us to the starting pair.
For $(x,y) \in \piop$ means $y \in x+1$, and multiplying through by $x^{-1}$ gives us 
$x^{-1}y \in x^{-1} + 1$, which says the pair $(x^{-1},x^{-1}y)$ is in $\piop$.
In other words, the pair $(x,y)$ gives the 
pair $(x^{-1},x^{-1}y)$.
Applying this formula again to $(x^{-1},x^{-1}y)$, we get
$((x^{-1})^{-1},(x^{-1})^{-1}x^{-1}y)$, which reduces to $(x,y)$.

So we can produce the block containing $(x,y)$ by alternately applying the two rules.
Reversal/negation gives $(-y,-x)$.
Then consistency gives $((-y)^{-1},(-y)^{-1}(-x)) = (-y^{-1},y^{-1}x)$.
Reversal/negation gives 
$(-y^{-1}x,y^{-1})$, and consistency gives 
$((-y^{-1}x)^{-1},(-y^{-1}x)^{-1}y^{-1}) =
(-x^{-1}y,-x^{-1})$ from that.
Then reversal/negation gives
$(x^{-1},x^{-1}y)$, and consistency gives
$(x,xx^{-1}y) = (x,y)$.

 To see that any block $B$ containing $(1,y)$ has at most three pairs, let $B$ be such a block.
 We have that any pair in a block can be obtained from any other by alternately applying reversal/negation and consistency.
 Consistency applied to $(1,y)$ yields $(1^{-1},1^{-1}y) = (1,y)$.
 Applying reversal/negation to it gives $(-y,-1)$, consistency gives $(-y^{-1},y^{-1})$,
 and reversal/negation leaves us at 
 $(-y^{-1},y^{-1})$.
\end{proof}

Returning to our observation that structures with many pairs in the $\piop$ relation tend to be hyperfields, the reason for this is that sums of three elements are likely to be all of $H$, giving us associativity since 
$(x+y)+z = H = x+(y+z)$.
Our previous paper \cite{hobbyjun2024} defined the term ``large sums'' differently, so we will not reuse it here with another meaning.
Hopefully the correct definition of ``large sums'' will eventually be uncovered.

The quantities $m$ and $k$ in Lemma \ref{Large sum lemma} are easily extracted from the relation $\piop$ of a structure.
We had that $m$ was the largest number so that 
for all $w \in H^\times$ we had $|(w+1)-\{0\}| \geq m$.
Since $\piop$ is only defined on $H^\times$, we have that $m$ is the least number of $1$s in rows of the incidence matrix of $\piop$.
We also had that $k$ was largest so that for all $w \in H^\times$ we had $|\{u \in H^\times \colon w \in (u+1)\}| \geq k$.
So $k$ is the least number of $1$s in columns of the incidence matrix of $\piop$.

\begin{mythm} \label{large gives hyperfield theorem}
Let $H^\times$ be a finite abelian group written multiplicatively, with distinguished element $-1$ of order $1$ or $2$.
Let $\piop$ be a binary relation on $H^\times$ which consists of a union of blocks, and let $m$ and $k$ be defined as in Lemma \ref{Large sum lemma}.
Also let $H_\pi$ be defined as in Definition \ref{H_pi definition}.  
Then if $m+k > |H^\times|$, the structure 
$H_\pi$ is a hyperfield.
\end{mythm}

\begin{proof}
Assume the hypotheses.
By construction, $H^\times$ is an abelian multiplicative group and 
$0 \cdot x = x \cdot 0 = 0$ for all $x \in H_\pi$.
We also have $0 + x = x + 0 = x$ for all 
$x \in H_\pi$ by construction.
Then to show addition is commutative, it suffices to do so on $H^\times$.
So let $x,y \in H^\times$ be given.
We have $x + y = y(y^{-1}x + 1)$, and claim this equals 
$y + x = x(x^{-1}y + 1)$.
So let $z \in y(y^{-1}x + 1)$ be given, which is equivalent to 
$y^{-1}z \in (y^{-1}x + 1)$ and to
$(y^{-1}x,y^{-1}z) \in \piop$.
Since $\piop$ in a union of blocks, it is closed under consistency, and we get
$((y^{-1}x)^{-1},(y^{-1}x)^{-1}y^{-1}z) = 
(x^{-1}y,x^{-1}z) \in \piop$
This gives us $x^{-1}z \in x^{-1}y + 1$, 
which implies $z \in x(x^{-1}+1)$, as desired.
All of the steps can be reversed, giving us commutativity.

Distributivity of addition comes from the definition.
Note that $z(x+y) = zx + zy$ holds if any of $x$, $y$ and $z$ is $0$, so we may assume 
$x,y,z \in H^\times$.
Then $z(x+y) = zy(y^{-1}x + 1) = zy(y^{-1}z^{-1}zx + 1) = 
zy((zy)^{-1}zx + 1)= zx + zy$.

We now have that every $x \in H^\times$ has a unique inverse, $(-1)x$.
For $x + (-1)x = 1\cdot x + (-1)\cdot x
= (1 - 1)x \ni 0$.
And if
$0 \in x + y = y(y^{-1}x + 1)$ we get 
$0 \in y^{-1}x + 1$, which implies 
$y^{-1}x = -1$ by our definition of addition.
Multiplying through by $(-1)y$ gives
$(-1)x = y$, as desired.

As noted in Lemma \ref{do not need reversibility lemma}, we do not need to show reversibility.
This leaves us to show that addition is associative.
If any of $x$, $y$ or $z$ is $0$, it is easy to show $(x+y)+z = x + (y+z)$, so we may assume $x,y,z \in H^\times$.
By Lemma \ref{Large sum lemma} we have that 
$H^\times \subseteq (x+y)+z$ and
$H^\times \subseteq x + (y+z)$, so it remains to show 
$0 \in (x+y)+z$ iff $0 \in x + (y+z)$.
Multiplying through by $y$ or $y^{-1}$, we see that it suffices to show 
$0 \in (x+1) + z$ iff $0 \in x+(1+z)$.

Assume $0 \in (x+1)+z$, so $-z \in (x+1)$ or
$x \piop -z$.
But $\piop$ is a union of blocks, and closed under reversal/negation.
Thus $z \piop -x$, which gives $-x \in z+1$
and $0 \in x + (z+1)$.
All of the steps are reversible, so addition is distributive.
\end{proof}

In view of Theorem \ref{large gives hyperfield theorem}, the condition that $m+k > r$ is important.
So we make the following definition.

\begin{mydef} \label{ample definition}
Let $H$ be a finite hyperfield, with $m$ and $k$ as in Lemma \ref{Large sum lemma}.
Then $H$ is {\em ample} iff $m + k > |H^\times|$.
\end{mydef}

\section{Numbers of Hyperfields}

In Section 5 of \cite{BakerJin2021}, Baker and Jin ask how the rate of growth of quotient hyperfields compares to that for all hyperfields.
We are now able to give a partial answer to this question.
It is convenient to use their notation, which defines $H_r$ to be the number of isomorphism classes of hyperfields of order $r+1$, and defines $Q_r$ to be the number of isomorphism classes of order $r+1$ that are quotient hyperfields.
They conjecture that the limit as $r$ goes to infinity of $Q_r/H_r$ is zero, based on a polynomial bound on the number of isomorphism classes of hyperfields of order $r+1$ that are quotients of finite fields.
We can prove this is true when $r$ is restricted to odd integers.

In outline, our argument will go as follows.
Given $r$, we consider all the possible blocks on finite abelian groups $H^\times$ of order $r$.
The structure of the blocks depends on the isomorphism class of $H^\times$.
By Theorem \ref{block size at most 6 Theorem}, the number $b$ of potential blocks is at least $r^2/6$.  Each of these may or may not be subsets of the relation $\piop$.
We can express which blocks are present with a vector $\vec{x} \in \{0,1\}^b$, where a component is $1$ if the corresponding block is in $\piop$.
By Theorem \ref{large gives hyperfield theorem} and the paragraph preceding it, $\piop$ gives a hyperfield structure on $H^\times$ whenever the number of $1$s in each of the rows and columns of its incidence matrix is greater than $r/2$.
These conditions on the numbers of $1$s give linear inequalities involving the components of $\vec{x}$, where each linear inequality is solved by around half of the vectors $\vec{x} \in \{0,1\}^b$.
Since there are at most $2r$ many such linear inequalities, we get that there are at least $2^b / 2^{2r}$ hyperfields of order $r+1$, where $b \geq r^2/6$.
This yields an exponential number of hyperfields, with a polynomial bound on how many are quotients of finite fields.

Unfortunately, we must also consider quotients of infinite fields.
As was proved by Bergelson and Shapiro in \cite{BS92} and also by Turnwald in \cite{T94}, if $G$ is a subgroup of finite index of the multiplicative group of an infinite field $F$, we have $G - G = F$.  This is also proved in \cite{BakerJin2021}.
This implies that $1-1 = H$ in every hyperfield $H$ that is a quotient of an infinite field.
Thus the only hyperfields of order $r+1$ that might be quotients of infinite fields are those where the row for $(-1)+1$ in the incidence matrix of $\piop$ consists of all $1$s.
This sometimes gives an upper bound of $2^{b-r}$ for the number of these quotient hyperfields, and that bound is quite close to our lower bound on $H_r$.

It is natural to label the rows and columns of the incidence matrix of $\piop$ by elements of $H^\times$.  
So we will write ``(-1) row'' to mean the row for $(-1)+1$, and so on.

So we will restrict most of our analysis to the best possible case, where $r$ is odd.
We state some lemmas for more general assumptions, in case that is useful later.
Otherwise, assume $r$ is odd in all of the following.

Since $r$ is odd, then -1 must have order 1 in the multiplicative group $H^\times$, giving $-1 = 1$.
Now every generator of $H^\times$ has order at least $3$, so $H^\times$ must have a generating set of size at most $\log_3(r)$.
Since an automorphism of a hyperfield is determined by the images of the generators of its multiplicative group, we have that isomorphism classes of hyperfields $H_\pi$ constructed as in Theorem \ref{large gives hyperfield theorem} have at most $r^{log_3(r)}$ elements.

Now let us investigate in more detail the conditions on which blocks must be present in $\piop$ in order that Theorem \ref{large gives hyperfield theorem} can be used to produce a hyperfield $H_\pi$.  Fix the group structure of $H^\times$, so we have a definite number $b$ of potential blocks, and also fix an ordering of these blocks, writing them as $B_1, B_2, B_3, \dots B_b$.
Then for a given $\vec{x} \in \{0,1\}^b$ we have an associated structure $H_\pi$, where $B_i \subseteq \piop$ if $\vec{x}(i) = 1$, and $B_i \cap \piop = \emptyset$ if $\vec{x}(i) = 0$.
Picking any $g \in H^\times$, the sum of the $g$ row in the incidence matrix of $\piop$ is equal to the dot-product 
$\vec{c}_g \cdot \vec{x}$, where $\vec{c}_g(i)$ is the number of pairs in $B_i$ that are in that row.
(We denote the $i$-th component of a vector $\vec{z}$ by 
$\vec{z}(i)$, since we will be using subscripts for other purposes.)
We now have the following definition.

\begin{mydef} \label{c_g definition}
Fix a group $H^\times$, and identify which element is $-1$.
Also fix a definite ordering of the blocks of $H^\times$, listing them as $B_1, B_2 \dots B_b$.
For $g \in H^\times$, define $\vec{c}_g$ to be the $b$-long vector where $\vec{c}_g(i)$ is the number of pairs in the $i$-th block of the form $(g,y)$.    
\end{mydef}

We do not need to use vectors for the columns of the incidence matrix of $\piop$, as they are the same vectors as for the rows.

\begin{lem} \label{columns same as rows lemma}
As before, fix a group $H^\times$, an element $-1$, and an ordering of the blocks.
Similarly to Definition \ref{c_g definition}, we define vectors $\vec{d}_g$ for the columns.  
For each $g \in H^\times$, define $\vec{d}_g$ to be the $b$-long vector where $\vec{d}_g(i)$ is the number of pairs in the $i$-th block of the form $(x,g)$.    
Then $\vec{d}_g = \vec{c}_{-g}$
\end{lem}

\begin{proof}
We have a one-to-one correspondence between pairs in the $g$ column and pairs in the $-g$ row, where $(x,g)$ corresponds to $(-g,-x)$.
The two pairs are related by reversal/negation, so they are in the same block.
This gives us that for any given block, the number of its pairs in the $g$ column and 
$-g$ row are equal, which implies 
$\vec{d}_g = \vec{c}_{-g}$.
\end{proof}

When $-1 = 1$, the situation is simpler. 
The reversal/negation operation is just reversal, making each $(x,g)$ in the same block as $(g,x)$. 

Since $r$ is odd and $m = k$ by 
Lemma \ref{columns same as rows lemma},
we will take both $m$ and $k$ to be $(r+1)/2$, giving $m+k > r$ in Theorem \ref{large gives hyperfield theorem}, as is required to show $H_\pi$ is a hyperfield.

So for the $g$ row of the incidence matrix of $\piop$ we have the inequality 
$\vec{c}_g \cdot \vec{x} > r/2$, and all of these row inequalities must be satisfied to guarantee that $H_\pi$ is a hyperfield.
The inequalities for columns are the same as for rows, and may be ignored.

\begin{lem} \label{two rows same inequality lemma}
For all $g \in H^\times$, $\vec{c}_g = \vec{c}_{g^{-1}}$.
\end{lem}

\begin{proof}
It suffices to show that each block $B_i$ has the same number of pairs in the $g$ row as in the $g^{-1}$ row.
The consistency operation on pairs exchanges $(x,y)$ and $(x^{-1},x^{-1}y)$
So each pair in $B_i$ in the $g$ row corresponds to a unique pair in the $g^{-1}$ row, and this one-to-one correspondence makes the number of pairs from $B_i$ the same in both rows.
\end{proof}

This is good when $r$ is odd, since it cuts the number of inequalities that $\vec{x}$ must satisfy roughly in half.
In general, we have the following lemma.

\begin{lem} \label{number of inequalities lemma}
Suppose $H^\times$ has $r$ elements, where $r = 2^kq$ with $q$ odd.  Then there are at most $(r+2^k)/2$ many distinct vectors $\vec{c}_g$ for $g \in H^\times$.
\end{lem}

\begin{proof}
Let $t$ be the number of elements of $H^\times$ of order $1$ or $2$, these elements are their own inverses.
Lemma \ref{two rows same inequality lemma} does not guarantee that these elements pair with other elements to give the same vector $\vec{c}_g$.
Thus the number of distinct $\vec{c}_g$ is bounded above by $t + (r-t)/2$.
By the Fundamental Theorem of Finite Abelian Groups, $H^\times$ is isomorphic to a product $T \times Q$, where $T$ is a group of order $2^k$ and $Q$ has odd order.
All elements of order $1$ or $2$ in this product have the form $(\tau,1_Q)$, giving us $t \leq 2^k$.
Thus the number of distinct $\vec{c}_g$ is at most $2^k + (r - 2^k)/2 = (r+2^k)/2$.
\end{proof}

Summarizing what we have so far when $r$ is odd, a vector $\vec{x} \in \{0,1\}^b$ gives us a hyperfield of even order $r+1$ whenever a set of inequalities of the form 
$\vec{c}_g \cdot \vec{x} > r/2$ is satisfied.
An irredundant list of the coefficient vectors $\vec{c}_g$ contains $(r+1)/2$ elements, there is $\vec{c}_1$ and $(r-1)/2$ other vectors $\vec{c}_g$ where one element $g$ is taken from each pair $\{h,h^{-1}\}$ from the non-identity elements of $H^\times$.

Before we proceed with the proof, an example is in order.  

\begin{myeg} \label{Z_7 example}
Let us take $r = 7$, so $H^\times$ is the $7$-element cyclic group generated by an element $a$.  In small examples, we label blocks with capital letters, and attempt to assign them in a canonical order.
\vspace{.5 cm}
\begin{center}
$\begin{array}{|c||c|c|c|c|c|c|c|} 
\hline
\piop & 1 & a & a^2 & a^3 &a^4 & a^5 & a^6\\
\hline
\hline
1 & A & B & C & D & E & F & G \\
\hline
a & B & G & H & I & J & K & H \\
\hline
a^2 & C & H & F & K & L & L & I \\
\hline
a^3 & D & I & K & E & J & L & J \\
\hline
a^4 & E & J & L & J & D & I & K \\
\hline
a^5 & F & K & L & L & I & C & H \\
\hline
a^6 & G & H & I & J & K & H & B \\
\hline
\end{array}$
\end{center}
\vspace{.5 cm}

The table shows the blocks, and is symmetric about its main diagonal.  
One sees that the $a$ row contains pairs from the blocks in $\{B, G, H, I, J, K\}$, with $H$ listed twice.
The $a^6$ row contains the same multiset of letters, but in a different order.
This continues, pairing the $a^2$ row and the $a^5$ row, and so on.

There are 12 different blocks.
Block A has size $1$, B through G each have size $3$, and blocks I through L each have size $6$.
This distribution of sizes is typical.

We stack the $4$ rows in our irredundant list of the $\vec{c}_g$, giving a $4 \times 12$ matrix $C$, and write $\vec{x}$ as a $12$-long column vector.  
Then the matrix form of our system of inequalities is 
$C \quad \vec{x} \quad > (3.5, 3.5, 3.5, 3.5)^T$, where each component of $C \vec{x}$ must be greater than $3.5$.
We of course let the first component of $\vec{x}$ denote the number of times block A appears, and continue in alphabetical order.
With these conventions, C is the following matrix.
$$
\begin{pmatrix}
1 & 1 & 1 & 1 & 1 & 1 & 1 & 0 & 0 & 0 & 0 & 0 \\
0 & 1 & 0 & 0 & 0 & 0 & 1 & 2 & 1 & 1 & 1 & 0 \\
0 & 0 & 1 & 0 & 0 & 1 & 0 & 1 & 1 & 0 & 1 & 2 \\
0 & 0 & 0 & 1 & 1 & 0 & 0 & 0 & 1 & 2 & 1 & 1 \\
\end{pmatrix}
$$

The sum of the entries in each row is the total number of times any block appears in the corresponding row of the incidence matrix for $\piop$, so it is of course our value of $r$, which is $7$.
\end{myeg}

For each row $\vec{c}_g$ of $C$, exactly half the vectors 
$\vec{x} \in \{0,1\}^b$ have dot product with $\vec{c}_g$ greater than $r/2$.
To see this, let $\hat{1}$ be the $b$-long vector of $1$s.
Then $\vec{c}_g \cdot \hat{1} = r$, and thus exactly one vector in each pair of $\vec{x}$ and $\hat{1} - \vec{x}$ has dot product with $\vec{c}_g$ that is greater than $r/2$, since
$\vec{c}_g \cdot \vec{x} + \vec{c}_g \cdot (\hat{1} - \vec{x}) = r$.

Since the angle in $\mathbb{R}^b$ between any two rows $\vec{c}_g$ and $\vec{c}_h$ is less than $90^\circ$, one would hope on geometric grounds that the proportion of vectors $\vec{x}$ with dot products greater than $r/2$ with all $d = (r+1)/2$ rows of $C$ would be at least $(1/2)^d$.
This is true, but there does not seem to be a nice proof of it.
Our argument will be algebraic, but could be interpreted geometrically as follows.  We first add extra dimensions to $\vec{x}$ and $C$, increasing the dimension from $b$ to $b'$, to ``give us room to work''.  Then we gradually ``rotate'' the hyperplanes about the center of the hypercube 
$\{0,1\}^{b'}$ until they are all mutually perpendicular, doing this in a way where we can demonstrate that each rotation decreases the number of solutions.
The number of solutions are easily counted when the hyperplanes are all perpendicular, giving us a lower bound on the original number of solutions.

\begin{mydef}
Let $C_1$ and $C_2$ be $(r+1)/2 \times b'$ matrices.  If all of the following hold, we say that $C_2$ is obtained from $C_1$ by a {\em valid $u$-$v$ swap on the $g$ row}.
\begin{enumerate}
  \item All of the entries of both matrices are nonnegative.
  \item The matrices are identical, except for two entries, those in the $g$ row and in the $u$ column and $v$ column.
  \item $C_1(g,u) = p$, $C_1(g,v) = 0$, $C_2(g,u) = 0$, and $C_2(g,v) = p$, for some positive real number $p$.
  \item The $v$ column of $C_1$ is all zeroes, and the $v$ column of $C_2$ is all zeroes except for a single $p$ in the $g$ row.
\end{enumerate}
\end{mydef}

This is the key lemma letting us do the ``rotations''.
We state it in more generality than is needed, since it may be of independent interest.

\begin{lem} \label{swap lemma}
Let $\vec{d}$ be a fixed column vector of length $(r+1)/2$.
Let $C_2$ be obtained from $C_1$ by a valid $u$-$v$ swap on the $g$ row.
Then the number of solutions with $\vec{x} \in \{0,1\}^{b'}$ of the system $C_1 \vec{x} > \vec{d}$ is at least as large as the number of solutions of the system $C_2 \vec{x} > \vec{d}$.
\end{lem}

\begin{proof}
As in \cite{Clark2003}, let us use the following notation for vectors obtained by minor modifications of the vector $\vec{x}$.
We write $\vec{x}^{\;u\,v}_{\;p\,q}$ for the vector that is equal to $\vec{x}$ at all components, except that its $u$-th component is $p$ and its $v$-th component is $q$.
We will also write $\vec{x} = \vec{x}^{\;u\,v}_{\;p\,q}$ as a compact way to say that $\vec{x}(u) = p$ and $\vec{x}(v) = q$.

Assume that $C_2$ is obtained from $C_1$ by a valid $u$-$v$ swap in the $g$ row, where $p$ is a positive real so that $C_1(g,u) = p$, $C_1(g,v) = 0$, $C_2(g,u) = 0$, and $C_2(g,v) = p$.  

To show that the system $C_1 \vec{x} > \vec{d}$ has at least as many solutions as $C_2 \vec{x} > \vec{d}$, we will produce an injective map $\phi$ from the set of solutions for $C_2$ to the set of solutions for $C_1$.
Let $\vec{x} \in \{0,1\}^{b'}$ be a given solution of 
$C_2 \vec{x} > \vec{d}$.

Our first case is where 
$\vec{x} \neq \vec{x}^{\;u\,v}_{\;0\,1}$.
We will show $\vec{x}$ is also a solution of 
$C_1 \vec{x} > \vec{d}$, so we may simply let $\phi(\vec{x})$ be $\vec{x}$.
We have $c_h \cdot \vec{x} > \vec{d}(h)$ for all rows $\vec{c}_h$ except $\vec{c}_g$, and we only need to worry about 
$\vec{c}_g \cdot \vec{x}$ in $C_1$ and $C_2$.
The change to $\vec{c}_g$ in going from $C_2$ to $C_1$ is to increase the $u$-th entry by $p$ and to decrease the $v$-th entry by $p$.  
The only way that this can decrease the dot product
$\vec{c}_g \cdot \vec{x}$ is if $\vec{x}(u) = 0$ and $\vec{x}(v) = 1$, which is not true in this case.

In the second case, assume we have $\vec{x}(u) = 0$ and $\vec{x}(v) = 1$, but that $\vec{c}_g \cdot \vec{x} > \vec{d}(g)+p$ in $C_2$.  
Then arguing as in the first case, we have that 
$\vec{x}$ is also a solution of $C_1 \vec{x} > \vec{d}$,
since the change in $\vec{c}_g \cdot \vec{x}$
upon going from $C_2$ to $C_1$ is only a decrease of $p$.
So again we let $\phi(\vec{x})$ be $\vec{x}$.

Our final case is where $\vec{x}(u) = 0$, $\vec{x}(v) = 1$, and $\vec{c}_g \cdot \vec{x} \leq \vec{d}(g)+p$ in $C_2$.
Here, we take $\phi(\vec{x})$ to be $\vec{x}^{\;u\,v}_{\;1\,0}$, the vector with the $u$-th and $v$-th components switched.
Since we have $c_h(v) = 0$ for all rows where $h \neq g$, we have $c_h \cdot \vec{x}^{\;u\,v}_{\;1\,0} \geq c_h \cdot \vec{x} > \vec{d}(h)$ in those rows.
In going from $C_2$ to $C_1$, the $u$-th and $v$-th components of $\vec{c}_g$ switch, and this is exactly the change made in going from $\vec{x}$ to $\vec{x}^{\;u\,v}_{\;1\,0}$.
Thus $\vec{c}_g \cdot \vec{x}^{\;u\,v}_{\;1\,0}$ in $C_1$ has the same value as $\vec{c}_g \cdot \vec{x}$ in $C_2$, and is also greater than $r/2$.
This shows that $\vec{x}^{\;u\,v}_{\;1\,0}$ is a solution of $C_1 \vec{x}^{\;u\,v}_{\;1\,0} > (r/2) \hat{1}$.

We also have that $\phi$ is injective, since in the third case the vector $\vec{x}^{\;u\,v}_{\;1\,0}$ is not a solution of $C_2 \vec{x}^{\;u\,v}_{\;1\,0} > \vec{d}$.
To see this, note that 
$\vec{c}_g \cdot \vec{x}^{\;u\,v}_{\;1\,0} = \vec{c}_g \cdot \vec{x} - p \leq \vec{d}(g)$ in $C_2$.
\end{proof}

Now we can prove the following.

\begin{mythm} \label{many hyperfields theorem}
Let $r$ be odd, and let $H^\times$ be an abelian group of order $r$, written multiplicatively.  Let $b$ be the number of blocks of the relation $\piop$ on $H^\times$.
Then at least $2^{b-(r+1)/2}$ of the structures $H_\pi$ are hyperfields.
\end{mythm}

\begin{proof}
Assume the hypotheses.
By Theorem \ref{large gives hyperfield theorem}, any relation $\piop$ that is a union of blocks making $m,k \geq (r+1)/2$, gives us $m + k \geq r+1 > r$ and makes the structure $H_\pi$ a hyperfield.
By Lemmas \ref{columns same as rows lemma} and \ref{two rows same inequality lemma}, we have such $m$ and $k$ whenever a set of $(r+1)/2$ inequalities holds for the vector $\vec{x}$ of coefficients of the possible blocks.
Specifically, each inequality is of the form 
$c_h \cdot \vec{x} > r/2$ where $c_h$ is the vector listing the number of times members of each possible block appear in the $h$ row of the incidence matrix.
We stack the $(r+1)/2$ distinct vectors $c_h$ to make the $(r+1)/2 \times b$ matrix $C$, where our inequalities combine to $C \vec{x} > \vec{d}$, with $\vec{d}(i) = r/2$ for all $i$.

Preparatory to a sequence of applications of Lemma \ref{swap lemma}, to the matrix $C$, we take $b'$ to be the total number of nonzero entries in $C$, and produce a matrix $C_1$
by adjoining $b'-b$ many columns of zeroes to the end of $C$.
This allows us to repeatedly swap entries in various rows until we arrive at a matrix where every column has at most one nonzero entry.

Specifically, we define $\nu(\vec{v})$ for any vector $\vec{v}$ to be one less than the number of nonzero elements in $\vec{v}$ if $\vec{v}$ has two or more nonzero elements, and define $\nu(\vec{v})$ to be zero otherwise.
So $\nu(\vec{v})$ is the number of ``extra'' nonzero elements in $\vec{v}$.
For any matrix $M$, we define $\gamma(M)$ to be the sum of $\nu(\vec{y})$ over all columns $\vec{y}$ of $M$.
So $\gamma(M)$ is the number of swaps needed to turn $M$ into a matrix with at most one nonzero entry in each column. 

Now we start with $C_1$.
If $\gamma(C_1) = 0$, we are done.
Otherwise, there is a column of $C_1$ with two or more nonzero entries, call it the $u$ column.
By construction, there is a column of all zeroes left in $C_1$, call it the $v$ column.
Pick a nonzero entry in the $u$ column, and say it is in the $g$ row of $C_1$.
Then there is a unique valid $u-v$ swap on the $g$ row of $C_1$ taking it to a matrix $C_2$.
All of the entries of $C_2$ are nonnegative, and $\gamma(C_2) = \gamma(C_1) - 1$.

We now do the same process to $C_2$, and continue until we reach the matrix $C_f$ where $\gamma(C_f) = 0$.
(In fact, $f = b' - b + 1$.)
By repeated applications of Lemma \ref{swap lemma}, we have that the number of solutions of $C_1 \vec{x} > \vec{d}$ for $\vec{x} \in \{0,1\}^{b'}$ is greater than or equal to the corresponding number of solutions of 
$C_f \vec{x} > \vec{d}$.

Valid $u-v$ swaps on any row simply rearrange elements of that row, so the sum of the entries in each row is still $r$.
Rearranging the columns of $C_f$ does not affect the number of solutions of $C_f \vec{x} > \vec{d}$, so we may rearrange the columns of $C_f$ to produce a block matrix 
$$D = \left(D_1 \vert D_2 \vert D_3 \vert \dots D_{(r+1)/2} \right)$$
Here, each matrix $D_i$ only has nonzero entries in the $i$-th row, all of the entries in that row are nonzero, and the sum of those entries is $r$.

A vector $\vec{x} \in \{0,1\}^{b'}$ is a solution of $D \vec{x} > \vec{d}$ iff it is a concatenation of vectors 
$\vec{x}_1, \vec{x}_2, \dots \vec{x}_{(r+1)/2}$
where each $\vec{x}_i$ is a solution of 
$R_i \vec{x}_i > r/2$ and each $R_i$ is the $i$-th row of $D_i$.
Consider any particular matrix $D_i$, and let $b_i$ be its number of columns.
Then exactly half of the vectors in $\{0,1\}^{b_i}$ are solutions of $R_i \vec{x}_i > r/2$.  
To see this, let $\hat{1}$ be the vector of length $b_i$ with all components $1$.
Each vector $\vec{y} \in \{0,1\}^{b_i}$ is paired with the vector $\hat{1} - \vec{y}$, where 
$R_i \vec{y} + R_i (\hat{1} - \vec{y}) = R_i \hat{1} = r$.
Thus exactly one element of each pair is a solution of 
$R_i \vec{x}_i > r/2$.
So for each $i$, $R_i \vec{x}_i > r/2$ has $2^{b_i -1}$ solutions, and 
$D \vec{x} > \vec{d}$ has 
$2^{b_1 -1} \cdot 2^{b_2 -1} \cdot \dots 2^{b_{(r+1)/2} -1} = 2^{b'- (r+1)/2}$ solutions. 

This implies that the original inequality $C_1 \vec{x} > \vec{d}$ has at least $2^{b'- (r+1)/2}$ solutions as well.
Recall that $C_1$ was produced from $C$ by adjoining $b'-b$ many columns of zeroes.
This means that there are $2^{b'-b}$ solutions of 
$C_1 \vec{x} > \vec{d}$ for every solution of
$C \vec{x} > \vec{d}$, giving us that there are at least 
$2^{b'- (r+1)/2}/2^{b'-b} = 2^{b-(r+1)/2}$ solutions of the latter inequality.
\end{proof}

A much easier inequality is the following.

\begin{mythm} \label{quotients of infinite fields theorem}
Fix a positive odd integer $r$ and an abelian group $H^\times$ of order $r$.
Since $r$ is odd, any hyperfield with multiplicative group $H^\times$ has $-1 = 1$.
Let $b$ be the number of possible blocks of $H^\times$.
Then there are at most $2^{b-r}$ hyperfields with the multiplicative group $H^\times$ that are quotients of infinite fields.
\end{mythm}

\begin{proof}
As was proved by Bergelson and Shapiro in \cite{BS92} or by Turnwald in \cite{T94}, if $G$ is a subgroup of finite index of the multiplicative group of an infinite field $F$, we have $G - G = F$.  This is also proved in \cite{BakerJin2021}.
This implies that $1-1 = H$ in every hyperfield $H$ that is a quotient of an infinite field.

Let $H$ be a hyperfield where $r$, the order of $H^\times$, is odd.  
Then $1 - 1 = 1 + 1$ in $H$, and to guarantee that $H$ is not a quotient of an infinite field it suffices to insure that some element of $H^\times$ is not in $1+1$.
By Theorem \ref{block size at most 6 Theorem}, each block of $H^\times$ contains at most one pair of the form $(1,y)$, and thus only hyperfields that contain all $r$ of these blocks could be quotients of infinite fields.  
Assuming that all the remaining blocks can be chosen arbitrarily, we get an upper bound of $2^{b-r}$ many possible hyperfields that are quotients of infinite fields.
\end{proof}

 In practice, many choices of sets of blocks do not yield hyperfields $H_\pi$, and some of those that do are isomorphic.

\begin{mythm} \label{limit of ratio theorem}
Let $H_r$ be the number of isomorphism classes of hyperfields of order $r+1$, and let $Q_r$ be the number of isomorphism classes of quotient hyperfields of order $r+1$.
Then when $r$ is restricted to odd integers, the limit as $r$ goes to infinity of $Q_r/H_r$ is zero.
\end{mythm}

\begin{proof}
Let $r$ be odd, and let $H^\times$ be any multiplicative abelian group of order $r$.
Since $r$ is odd, we have $-1 = 1$.
The exact number of blocks $b$ depends on the group $H^\times$, but we have $r^2/6 \leq b \leq r^2$ by Theorem \ref{block size at most 6 Theorem}.

By Theorem \ref{many hyperfields theorem}, we have that there are at least $2^{b-(r+1)/2}$ hyperfields $H_\pi$ with multiplicative group $H^\times$.
By the discussion at the start of this section, at most $r^{\log_3(r)} < 2^{\log_2(r)^2}$ of the $H_\pi$ can be in the same isomorphism class.
This gives us that there are at least 
$2^{b-(r+1)/2 - \log_2(r)^2}$ distinct isomorphism classes of hyperfields of order $r+1$.

By Theorem 1.1 of \cite{BakerJin2021}, when $r$ is odd an upper bound on the number of isomorphism classes of hyperfields that are quotients of finite fields is given by the number of fields of size less than or equal to $r^4$.
Although better bounds are possible, we will simply use $r^4$.

By our Theorem \ref{quotients of infinite fields theorem}, 
the number of isomorphism classes of hyperfields that are quotients of infinite fields is bounded above by $2^{b-r}$.
In the limit, we have $2^{b-r} \geq 2^{r^2/6-r} \geq r^4$, showing that the number of isomorphism classes of quotients of finite fields is also below 
$2^{b-r}$.
Thus the total number of isomorphism classes of quotient hyperfields is at most $2\cdot2^{b-r} = 2^{b-r+1}$.

Thus $Q_r/H_r \leq 2^{b-r+1} / 2^{b-(r+1)/2 - \log_2(r)^2} =
2^{-r + 1 +(r+1)/2 + \log_2(r)^2} = 2^{-r/2 + 3/2 + \log_2(r)^2}$, which goes to zero as $r$ goes to infinity.
\end{proof}

It is natural to ask if this result can be strengthened.
First note that factors smaller than $2^r$ are unlikely to matter in the limit.
Thus we can basically ignore the $r^4$ bound on the number of quotients of finite fields, and the $r^{\log(r)}$ bound on the size of isomorphism classes of the $H_\pi$.
Lemma \ref{number of inequalities lemma} writes $r = 2^kq$ with $q$ odd, and bounds the number of inequalities that need to be solved to produce hyperfields by $(r+2^k)/2$.
For any fixed $k$ this only introduces a constant factor, which is negligible.  
Each additional inequality seems to reduce the number of hyperfields by a factor of $2$, so it appears at first that the theorem might still hold for small values of $k$.

Unfortunately, even when $r$ has a single factor of $2$ we no longer have that $-1$ must be $1$.
In Theorem \ref{quotients of infinite fields theorem}, we replaced $1 + (-1)$ with $1+1$, and used Theorem \ref{block size at most 6 Theorem} to argue that $r$ different blocks needed to be included in $H_\pi$ in order to make $1+1 = H$.
This gave us the upper bound of $2^{b-r}$ on the number of hyperfields that were quotients of infinite fields.
It is no longer true that pairs of the form $(-1,y)$ take up $r$ blocks when $-1 \neq 1$.
For example, given any pair $(-1,y)$, consistency implies that $((-1)^{-1},(-1)^{-1}y) = (-1,-y)$ is in the same block.
This reduces the number of blocks needed to make $(-1)+1$ to $r/2$, giving us an upper bound of $2^{b-r/2}$ on the number of hyperfields that are quotients of infinite fields.
This is too large for the proof to work.

We could salvage a little by restricting to hyperfields with $-1 = 1$, writing $r = 2^kq$ and keeping $k < \log_2(r)$.
But it would be better to have another tool to exclude quotients of infinite fields.

\section{FETVINS for ample hyperfields}

In \cite{baker2022some}, M.~Baker and T.~Zhang asked whether all hyperfields had a certain property.
In our earlier paper \cite{hobbyjun2024}, we called this property FETVINS. 

\begin{mydef} \label{FETVINS definition}
We say that a hyperfield $F$ has the 
{\em FETVINS} property iff 
``Fewer Equations Than Variables Implies Nontrivial Solutions''.
That is, whenever we have a system of $k$ homogeneous linear equations in $n$ variables where $k < n$, there are values for the variables $x_1, x_2, \dots x_n$ which are not all zero, so that all $k$ of the homogeneous linear equations are true.
\end{mydef}

It is known that all quotient hyperfields have FETVINS, they inherit the property from the fields they are quotients of.
In \cite{hobbyjun2024}, we showed that a class of hyperfields that contained many nonquotient hyperfields also had FETVINS.
This class had what we called ``large sums'', a property similar to what we called ample in Definition \ref{ample definition}.
Thus it is natural to ask if all ample hyperfields have FETVINS, particularly as many finite hyperfields are ample.

The answer is positive, and given by the following theorem.

\begin{mythm} \label{ample implies FETVINS theorem}
Every finite hyperfield that is ample has FETVINS.
\end{mythm}

The proof of this theorem is similar to that of our Theorem 3.11 in \cite{hobbyjun2024}.
A key lemma is the following, it lets us ignore equations with four or more terms.

\begin{lem} \label{sum of 4 lemma}
In a finite hyperfield that is ample, any sum of four or more nonzero elements is the entire hyperfield.
\end{lem}

\begin{proof}
Let $H$ be a finite hyperfield that is ample.
As is well known, $H + x = H$ for any $x \in H$.
(Given $y \in H$, we have $y = y+0 \subseteq y + (-x + x)
= (y-x) +x \subseteq H+x$.)
So it is enough to show that the sum of any four nonzero elements is all of $H$.

By Lemma \ref{Large sum lemma}, the sum of any three elements of $H^\times$ contains $H^\times$.
Given $a,b,c,d \in H^\times$, we have 
$(a+b+c)+d \supseteq H^\times + d$, and claim that 
$H^\times + d = H$.
We have $-d \in H^\times$, so $0 \in H^\times + d$.
Now let $e \in H^\times$ be given.
Then $d^{-1}e \in H^\times$, and we have 
$d^{-1}e \in 1+f$ for some $f \in H^\times$ since $H$ is ample.
Thus $e \in d + df \subseteq H^\times + d$, showing 
$H^\times + d = H$.
\end{proof}

Another lemma helps us with equations with three terms.

\begin{lem} \label{solve two at once lemma}
Let $H$ be a finite hyperfield that is ample, and
let $a,b,c,d \in H^\times$ be given.
Then there exists $e \in H^\times$ so that both
$0 \in a+b+e$ and $0 \in c+d+e$.
\end{lem}

\begin{proof}
Let $H$ be ample, with $a,b,c,d \in H^\times$.
Although we originally defined $m$ and $k$ separately, Lemma \ref{columns same as rows lemma} gives us that $m$ and $k$ are equal.
For $m$ is the minimum value of the sum of the components of vectors of the form $\vec{c}_g$, while $k$ is the minimum value of the sum of components of column vectors  $\vec{d}_h$.
Since $H$ is ample, $m+k > |H^\times|$, or $m+m > |H^\times|$.
Now $a+b = a(1+a^{-1}b)$ has at least $m$ nonzero elements.
And similarly, $c+d$ has at least $m$ nonzero elements.
Thus $a+b$ and $c+d$ have some nonzero element $f$ in common.
We take $e$ to be $-f$, giving $a+b+e \supseteq f-f \ni 0$ and $c+d+e \supseteq f-f \ni 0$.
\end{proof}

Our proof of Theorem \ref{ample implies FETVINS theorem} now proceeds as follows.
Let a system of the following form be given, where $m < n$.

\begin{align*}
a_{11} x_1 + a_{12} x_2 + \dots a_{1n} x_n &\ni 0\\
a_{21} x_1 + a_{22} x_2 + \dots a_{2n} x_n &\ni 0\\
\vdots\\
a_{m1} x_1 + a_{m2} x_2 + \dots a_{mn} x_n &\ni 0
\end{align*}

We call the statements in this system {\em homogeneous linear equations}, since that is what they are analogous to.
Suppose that this system is a counterexample to FETVINS in $H$ that has the least number of variables.
That is, it has fewer equations than variables and no nontrivial solutions.

Some of the coefficients $a_{ij}$ may be zero.
We adopt the convention that only terms with nonzero coefficients are considered.
If one of our equations has only a single nonzero coefficient, we say that it has a {\em single variable} and write it in the form $a_{ij} x_j \ni 0$.
No single variable equation can occur in our minimal counterexample, for we could simply set $x_j$ equal to zero in our system, removing that variable and equation.
This would give us a smaller counterexample, since any nontrivial solution of the smaller system gives us a nontrivial solution of the larger one.

We can also not have any two-variable equations.  
For suppose one of the equations is 
$a_{ij}x_j + a_{ik} x_k \ni 0$.
Then we substitute $-a_{ij}^{-1} a_{ik} x_k$ for all occurrences of $x_j$ in our system.
This solves the two variable equation, and removes the variable $x_j$, producing a smaller system. 
One possible issue with doing this is that if say equation $p$ contains $a_{pj}x_j + a_{pk} x_k$, then doing the substitution produces 
$a_{pj}(-a_{ij}^{-1} a_{ik}) x_k + a_{pk} x_k = 
(-a_{pj}(-a_{ij}^{-1} a_{ik}) + a_{pk}) x_k$,
where the coefficient of $x_k$ is not a single element.
If this happens, we can fix it by replacing 
$(-a_{pj}(-a_{ij}^{-1} a_{ik}) + a_{pk})$ with any of its nonzero elements, for then any nontrivial solution of the smaller system produces a nontrivial solution of the original system.

We need another definition.

\begin{mydef}
A {\em subsystem} of a system of homogeneous linear equations consists of some subset of those equations, together with the set of all variables occurring in those equations.  Call a subsystem a {\em pile} if it has at least as many equations as variables.
\end{mydef}

Then our minimal counterexample to FETVINS can not contain any piles, for if it did we would produce a smaller system by setting all of the variables in the pile equal to zero.
This would remove at least as many equations as variables, and result in a smaller counterexample to FETVINS.
(In fact, a one-variable equation is a pile, so this subsumes the case of one-variable equations.)

We have now established that our minimal counterexample has no one-variable or two-variable equations, and also has no piles.
We will show that this situation is impossible, by constructing a solution of the system where none of the variables is zero.
It is enough to do this for the subsystem consisting of the three-variable equations and their variables, since by Lemma \ref{sum of 4 lemma} any homogeneous equation with more than three variables is true for all nonzero values of its variables.

Let us say that our subsystem has $e$ many three-variable equations, and $v$ many variables.
We have $e < v$, so the number of occurrences of variables in our equations is $3e$, which is less than $3v$.
Thus some variable occurs $2$ or fewer times in the subsystem, without loss of generality call this variable $x_1$.
We now remove the one or two equations containing $x_1$, and still have more variables than equations.
So we find a variable which we may call $x_2$ that is in $2$ or fewer of the remaining equations, remove those equations, and repeat this all again on the smaller system.
At the end of this process, we have a sequence of variables which we may as well call $x_1, x_2, x_3 \dots x_p$, where each variable is in two or fewer equations which use just it and subsequent variables.

Now we construct our solution by working backwards from the ``$p$'' end of the sequence.
We start by letting $x_{p-1}$ and $x_p$ have arbitrary nonzero values.
There are now two or fewer equations involving only variables in the set $\{x_{p-2}, x_{p-1}, x_p\}$, so Lemma \ref{solve two at once lemma} lets us find a nonzero value for $x_{p-2}$ that gives a solution to all of them.
We continue, and since $x_{p-3}$ is involved in two or fewer equations involving variables that have been assigned values, another application of Lemma \ref{solve two at once lemma}, solves all of those equations.
Repeating the process, we produce the desired solution with all variables nonzero.

This completes our proof of Theorem \ref{ample implies FETVINS theorem}.

\bibliography{rank}\bibliographystyle{alpha}
\end{document}